\documentclass{amsart}
 
\RequirePackage{amsmath}
\RequirePackage{amssymb}
 
\theoremstyle{plain}
\newtheorem{theorem}{Theorem}[section]
\newtheorem{proposition}[theorem]{Proposition}
\newtheorem{lemma}[theorem]{Lemma}

\newtheorem{mainthm}{Theorem}

\theoremstyle{definition}
\newtheorem{definition}[theorem]{Definition}
\newtheorem{example}[theorem]{Example}
\newtheorem{remark}[theorem]{Remark}

\numberwithin{equation}{section}

\newcommand{\C}{\mathbf{C}}

\newcommand{\F}{\mathbf{F}}
\newcommand{\Fl}{\F_{\ell}}
\newcommand{\Flbar}{\bar{\F}_{\ell}}
\newcommand{\Kbar}{\bar{K}}
\newcommand{\Q}{\mathbf{Q}}
\newcommand{\Qbar}{\bar{\Q}}
\newcommand{\Ql}{\Q_{\ell}}
\newcommand{\Qlbar}{\bar{\Q}_{\ell}}
\newcommand{\Qp}{\Q_{p}}
\newcommand{\Qpbar}{\Qbar_{p}}

\newcommand{\GQ}{G_{\Q}}
\newcommand{\GQl}{G_{\Q,\{\ell\}}}
\newcommand{\GQS}{G_{\Q,S}}
\newcommand{\GQSl}{G_{\Q,S\cup\{\ell\}}}

\newcommand{\CC}{\mathcal{C}}
\newcommand{\m}{\mathfrak{m}}
\renewcommand{\O}{\mathcal{O}}
\newcommand{\Z}{\mathbf{Z}}
\newcommand{\Zl}{\Z_{\ell}}
\newcommand{\Zp}{\Z_{p}}

\newcommand{\G}{\mathcal{G}}
\newcommand{\Hf}{H^{1}_{f}}

\newcommand{\MF}{\mathcal{MF}}
\newcommand{\U}{\mathcal{U}}

\newfont{\cyrr}{wncyr10}
\newcommand{\Sha}{\mbox{\cyrr Sh}}

\newcommand{\Fbar}{\bar{F}}
\newcommand{\Klm}{K_{\lambda}}
\newcommand{\Klmbar}{\Kbar_{\lambda}}

\newcommand{\ilm}{\iota_{\lambda}}
\newcommand{\klm}{k_{\lambda}}
\newcommand{\klmbar}{\bar{k}_{\lambda}}
\newcommand{\Olm}{\O_{\lambda}}

\newcommand{\rhob}{\bar{\rho}}
\newcommand{\rhobfl}{\bar{\rho}_{f,\ell}}
\newcommand{\rhobflm}{\rhob_{f,\lambda}}

\newcommand{\rhoblm}{\rhob_{\lambda}}
\newcommand{\rhofl}{\rho_{f,\ell}}
\newcommand{\rhoflm}{\rho_{f,\lambda}}

\newcommand{\rholm}{\rho_{\lambda}}

\newcommand{\chib}{\bar{\chi}}

\newcommand{\chiblm}{\chib_{\lambda}}

\newcommand{\chilm}{\chi_{\lambda}}

\newcommand{\eplb}{\bar{\varepsilon}_{\ell}}
\newcommand{\eplm}{\varepsilon_{\lambda}}
\newcommand{\eplmb}{\bar{\varepsilon}_{\lambda}}
\newcommand{\nm}{\left| \cdot \right|}

\newcommand{\inj}{\hookrightarrow}

\newcommand{\surj}{\twoheadrightarrow}
\newcommand{\toi}{\overset{\simeq}{\longrightarrow}}
\newcommand{\too}{\longrightarrow}
\newcommand{\ts}{\textstyle}

\DeclareMathOperator{\ab}{ab}
\DeclareMathOperator{\ad}{ad}
\DeclareMathOperator{\cc}{c}
\DeclareMathOperator{\crys}{crys}
\DeclareMathOperator{\End}{End}
\DeclareMathOperator{\Frob}{Frob}
\DeclareMathOperator{\Gal}{Gal}
\DeclareMathOperator{\GL}{GL}
\DeclareMathOperator{\Hom}{Hom}

\DeclareMathOperator{\Ind}{Ind}
\DeclareMathOperator{\PGL}{PGL}
\DeclareMathOperator{\proj}{proj}
\DeclareMathOperator{\rank}{rank}

\DeclareMathOperator{\Sets}{Sets}
\DeclareMathOperator{\sss}{ss}

\DeclareMathOperator{\trace}{tr}

\newcommand{\adz}{\ad^{0}\!}

\begin{document}

\title{Unobstructed modular deformation problems}
\author{Tom Weston}
\email{weston@math.berkeley.edu}
\address{Department of Mathematics, University of California, Berkeley}
\thanks{Supported by an NSF postdoctoral fellowship}

\begin{abstract}
Let $f$ be a newform of weight $k \geq 3$ with
Fourier coefficients in a number
field $K$. We show that the universal deformation ring of the mod $\lambda$ 
Galois
representation associated to $f$ is unobstructed, and 
thus isomorphic to a power
series ring in three variables over the Witt vectors, for all but finitely many
primes $\lambda$ of $K$. 
We give an explicit bound on such $\lambda$ for the 6 known
cusp forms of level 1, trivial character, and rational Fourier coefficients. We
also prove a somewhat weaker result for weight 2.
\end{abstract}

\maketitle

\section{Introduction} \label{intro}

Let $f = \sum a_{n}q^{n}$ be a newform of weight $k \geq 2$, level
$N$, and character $\omega$; let $S$ be any finite set of places of $\Q$ 
containing the infinite place and all primes dividing $N$.
For every prime $\lambda$ of the number field $K$ generated by the $a_{n}$,
Deligne and Serre
have associated to $f$ a semisimple two dimensional Galois representation
$$\rhobflm : \GQSl \to \GL_{2} \klm$$
over the residue field $\klm$ of $\lambda$; here
$\GQSl$ is the Galois group of the maximal extension of $\Q$
unramified outside $S$ and the characteristic $\ell$ of $\klm$.
The representation
$\rhobflm$ is absolutely irreducible for almost all $\lambda$;
for such $\lambda$ let $R^{S}_{f,\lambda}$ denote the
universal deformation ring parameterizing lifts of $\rhobflm$
(up to strict equivalence) to two dimensional representations of
$\GQSl$ over
noetherian local rings with residue field $\klm$ (see Section~\ref{s31} for a
precise definition).  Using recent work of Diamond, Flach, and Guo
\cite{DFG} on the Bloch--Kato conjectures for adjoint motives of modular
forms we prove the following theorem.

\begin{mainthm} \label{mt1}
If $k > 2$, then 
\begin{equation} \label{eq:main}
R^{S}_{f,\lambda} \cong W(\klm)[[T_{1},T_{2},T_{3}]]
\end{equation}
for all but finitely many primes $\lambda$ of $K$ (depending on $f$ and $S$); 
here $W(\klm)$ is
the ring of Witt vectors of $\klm$.  If $k=2$, then
(\ref{eq:main}) holds for all but finitely many primes $\lambda$ of $K$
dividing rational primes $\ell$ such that
\begin{equation} \label{eq:cond}
a_{\ell}^{2} \not\equiv \omega(\ell) \pmod{\lambda}.
\end{equation}
\end{mainthm}

The special case of Theorem~\ref{mt1} for elliptic curves was proven by
Mazur \cite{Mazur2} 
using results of Flach \cite{Flach2} on symmetric
square Selmer groups of elliptic curves.  For weight $k \geq 3$,
Theorem~\ref{mt1} answers Mazur's
question of \cite[Section 11]{Mazur2} concerning the finiteness of
the set of obstructed primes for modular deformation problems; 
previously this finiteness
was not known for a single modular form.  We refer
to \cite{Mazur2} for a discussion of additional applications of
Theorem~\ref{mt1}.

Our methods are in principle effective: that is, given enough information
about the modular form $f$ it is possible to determine a finite set of primes
$\lambda$ containing all those violating (\ref{eq:main}).  As a first example,
we combine Theorem~\ref{mt1} with work of Hida and Mazur to prove the
following theorem.

\begin{mainthm} \label{mt2}
Let $f$ be one of  the normalized
cusp forms of level 1, weight $k=12,16,18,20,22\text{~or~}26$,
and trivial character.  Let $\ell > k+1$ be a prime for which
$\rhobfl$ is absolutely irreducible.  Then
\begin{equation} \label{eq:two}
R_{f,\ell}^{\emptyset} \cong \Zl[[T_{1},T_{2},T_{3}]].
\end{equation}
\end{mainthm}

For example, for $f = \Delta$, (\ref{eq:two}) holds for $\ell \geq 17$,
$\ell \neq 691$.

By \cite{Mazur}, to prove Theorem~\ref{mt1} it suffices to show that the
Galois cohomology group
$H^{2}(\GQ,\ad \rhobflm)$ of the adjoint representation of $\rhobflm$
vanishes for almost all $\lambda$ (or for almost all $\lambda$ satisfying
(\ref{eq:cond}) when $k=2$).  In Section~\ref{s3} we explain how to
use Poitou--Tate duality and results on Selmer groups as in \cite{DFG}
to reduce the proof of Theorem~\ref{mt1} to two statements:
\begin{enumerate}
\item For fixed $p$, $H^{0}(\Gal(\Qpbar/\Qp),\eplmb \otimes \ad \rhobflm) = 0$
for almost all $\lambda$;
\item $H^{0}(\Gal(\Qlbar/\Ql),\eplmb \otimes \ad \rhobflm) = 0$
for almost all $\lambda$ (or for almost all $\lambda$ satisfying
(\ref{eq:cond}) when $k=2$), where $\ell = \text{char~} \klm$.
\end{enumerate}

In the elliptic curve case, the restriction of $\rhob_{f,\ell}$ to the
absolute Galois group $G_{p}$ of $\Qp$ has a canonical interpretation
as the Galois action on the $\ell$-torsion points of the elliptic curve;
the proofs of (1) and (2) are then a straightforward exercise
using formal groups and the Kodaira--N\'eron classification of special
fibers of
N\'eron models of elliptic curves.  Unfortunately, for more general modular
forms there is
no canonical model of the (often reducible) representation
$\rhobflm|_{G_{p}}$, and the tools from the elliptic curve case are no
longer applicable.

Our proof of (1) divides into two cases. 
Let $\pi_{p}$ be the
$p$-component of the automorphic representation corresponding to $f$.
If $\pi_{p}$ is principal series or
supercuspidal, then there is no harm in studying the
semisimplification of $\rhobflm|_{G_{p}}$.  One can then use the local
Langlands correspondence for $\GL_{2}$ to deduce (1) for such $p$.
This is done in Section~\ref{s1}.
If $\pi_{p}$ is special, then to prove (1) for $p$ 
it is not sufficient to consider the semisimplification, and the lack of a
canonical choice for $\rhobflm|_{G_{p}}$ prevents one from using a purely
local approach.
Instead, we use a level-lowering argument suggested to us by Ken
Ribet to verify (1) for such $p$.  
We present this argument in Section~\ref{s41}.

We use the theory of Fontaine--Laffaille to reduce (2) to
a computation with filtered 
Dieudonn\'e modules; this is presented in Section~\ref{s2}.  
The proofs of Theorems~\ref{mt1} and~\ref{mt2} are given in Section~\ref{s4}.

It is a pleasure to thank Matthias Flach,
Elena Mantovan, Robert Pollack, and Ken Ribet for
helpful conversations related to this work.

\subsection*{Notation}

If $\rho : G \to \GL_{2} K$ is a representation of a group $G$ over a 
field $K$, we write $\ad \rho : G \to \GL_{4} K$ for the adjoint representation
of $G$ on $\End(\rho)$ and $\adz \rho : G \to \GL_{3} K$
for the kernel of the trace map from $\ad \rho$ to the trivial representation.

We write $\GQ := \Gal(\Qbar/\Q)$ for the absolute Galois group of $\Q$.
We fix now and forever embeddings $\Qbar \inj \Qpbar$ for each $p$; these
yield injections $G_{p} \inj \GQ$ with $G_{p} := \Gal(\Qpbar/\Qp)$
the absolute Galois group of $\Qp$.  All Frobenius elements are geometric,
although for a character $\omega$ unramified at $\ell$ 
we use the usual normalization
$\omega(\ell) := \omega(\Frob_{\ell}^{-1})$.
We use the
phrase ``almost all'' as a synonym for the phrase ``all but finitely many''.

\section{Deformation theory} \label{s3}

\subsection{Universal deformation rings} \label{s31}

Let $k$ be a finite field of odd characteristic $\ell$.
Let $S$ be a finite set of places of $\Q$ containing $\ell$ and 
the infinite place, and
let $\Q_{S}$ denote the maximal extension of $\Q$ unramified away from $S$;
set $\GQS := \Gal(\Q_{S}/\Q)$.
Consider an absolutely irreducible continuous Galois representation
$$\rhob : \GQS \to \GL_{2}k.$$
We further assume that $\rhob$ is {\it odd} in the sense that the
image of complex conjugation under $\rhob$ has distinct eigenvalues.

Let $\CC$ denote the category of inverse limits of artinian local rings
with residue field $k$; morphisms in $\CC$ are assumed to induce the identity
map on $k$.  For a ring $A$ in $\CC$, we say that
$\rho : \GQS \to \GL_{2}A$
is a {\it lifting} of $\rhob$ to $A$ if the composition
$$\GQS \overset{\rho}{\too} \GL_{2}A \to \GL_{2}k$$
equals $\rhob$.  We say that two liftings $\rho_{1}, \rho_{2}$ of $\rhob$
to $A$ are {\it strictly conjugate} if there is a matrix $M$ in the kernel
of $\GL_{2}A \to \GL_{2}k$ such that
$\rho_{1} = M\cdot\rho_{2}\cdot M^{-1}$.
A {\it deformation} of $\rhob$ to $A$ is a strict conjugacy class
of liftings of $\rhob$ to $A$.  The {\it deformation functor}
$$D^{S}_{\rhob} : \CC \to \Sets$$
sends a ring $A$ to the set of deformations of $\rhob$ to $A$.
Since $\rhob$ is absolutely irreducible, by \cite[Proposition 1]{Mazur}
the functor $D^{S}_{\rhob}$ is
represented by a ring $R_{\rhob}$ of $\CC$; that is, there is an
isomorphism of functors
$$D^{S}_{\rhob}(-) \cong \Hom_{\CC}(R_{\rhob},-).$$
We call $R_{\rhob}$ the {\it universal deformation ring} for the deformation
problem $D^{S}_{\rhob}$.

The deformation problem $D^{S}_{\rhob}$ is said to be {\it unobstructed} if
$H^{2}(\GQS,\ad \rhob)=0$.  In this case,
\cite[Sections 1.7 and 1.10]{Mazur} gives the following description of
$R_{\rhob}$.

\begin{proposition} \label{prop:defthy}
If $D^{S}_{\rhob}$ is unobstructed, then
$$R_{\rhob} \cong W(k)[[T_{1},T_{2},T_{3}]]$$
with $W(k)$ the ring of Witt vectors of $k$.
\end{proposition}

\subsection{A criterion for unobstructedness} \label{s32}

We continue with the notation of the previous
section.  In this section we recall Flach's criterion of
\cite[Section 3]{Flach2} for the vanishing of
$H^{2}(\GQS,\ad \rhob)$.
Let $\O \in \CC$ be a totally ramified integral extension of $W(k)$.  Let
$\m$ denote the maximal ideal of $\O$ and fix a lifting
$$\rho : \GQS \to \GL_{2}\O$$
of $\rhob$ to $\O$.  
Let $K$ denote the fraction field of $\O$ and
let $V_{\rho}$ (resp.\ $A_{\rho}$) denote a three dimensional $K$-vector
space (resp.\ $(K/\O)^{3}$) endowed with a $\GQS$ action
via $\adz \rho$; we write $V_{\rho}(1)$ and $A_{\rho}(1)$ for their
Tate twists.

We recall that the {\it Selmer group} of any these Galois modules $M$ is
defined by
$$\Hf(\GQ,M) := \bigl\{ c \in H^{1}(\GQ,M) \,;\, c|_{G_{p}} \in \Hf(G_{p},M)
\text{~for all~} p \bigr\}$$
where $\Hf(G_{p},M)$ is as in \cite[Section 3]{BlochKato}.
Let $V$, $A$ denote either of the pairs $V_{\rho}$, $A_{\rho}$ or
$V_{\rho}(1)$, $A_{\rho}(1)$.  We define $\Sha(A)$ via the exact sequence
\begin{equation} \label{eq:finite}
0 \to i_{*}\Hf(\GQ,V) \to \Hf(\GQ,A) \to \Sha(A) \to 0
\end{equation}
with $i : V \surj A$ the natural map.
Finally, for any $\GQS$-module $M$ we define
$$\Sha^{1}(\GQS,M) := \{ c \in H^{1}(\GQS,M) \,;\, c|_{G_{p}} = 0
\text{~for all~} p \in S\}.$$
Let $\eplb : \GQ \to k^{\times}$ denote the mod $\ell$ Teichm\"uller character.

\begin{proposition} \label{prop:criterion}
Let $\rhob$ be as above.  Suppose that
\begin{enumerate}
\item \label{pc:hyp1} 
$H^{0}(G_{p},\eplb \otimes \ad \rhob) = 0$ for all $p \in S - \{\infty\}$;
\item \label{pc:hyp2}
$\Hf(\GQ,V_{\rho}(1)) = 0$;
\item \label{pc:hyp3}
$\Hf(\GQ,A_{\rho})=0$.
\end{enumerate}
Then the deformation problem $D^{S}_{\rhob}$ is unobstructed.
\end{proposition}
\begin{proof}
The trace pairing on $\ad \rhob$ identifies $\eplb \otimes \ad \rhob$ with
the Cartier dual of $\ad \rhob$.
Poitou--Tate duality (see
\cite[Proposition 4.10]{Milne}) thus yields an exact sequence
\begin{multline} \label{eq:es}
\underset{p \in S - \{\infty\}}{\ts \prod}
H^{0}(G_{p},\eplb \otimes \ad \rhob) \to
\Hom_{k}\bigl(H^{2}(\GQS,\ad \rhob),k\bigr) \to \\
\Sha^{1}(\GQS,\eplb \otimes \ad \rhob).
\end{multline}
The latter group decomposes as
$$\Sha^{1}(\GQS,\eplb \otimes \ad \rhob) =
\Sha^{1}(\GQS,\eplb) \oplus \Sha^{1}(\GQS,\eplb \otimes \adz \rhob),$$
and the first summand vanishes by \cite[Lemma 10.6]{Weston}.  Thus,
by (\ref{eq:es}) and hypothesis (\ref{pc:hyp1}), to prove the proposition it
suffices to show that $\Sha^{1}(\GQS,\eplb \otimes \adz \rhob)=0$.

Since $\rho$ is a lifting of $\rhob$,
the $\GQS$-module $A_{\rho}(1)[\m]$ is a realization of
$\eplb \otimes \adz \rhob$.
Since $\rhob$ is irreducible, the natural map
$$H^{1}(\GQ,\eplb \otimes \adz \rhob) 
= H^{1}(\GQ,A_{\rho}(1)[\m]) \to H^{1}(\GQ,A_{\rho}(1))$$
is injective.
The image of $\Sha^{1}(\GQS,\eplb \otimes \adz \rhob)$
under this map is easily seen to lie in $\Hf(\GQ,A_{\rho}(1))$, so that
to complete the proof it suffices to show that this Selmer group vanishes.

By (\ref{eq:finite}) and hypothesis (\ref{pc:hyp3}) we have $\Sha(A_{\rho})=0$.
Since $A_{\rho}(1)$ is Cartier dual to $A_{\rho}$, by
\cite[Theorem 1]{Flach1} this implies that $\Sha(A_{\rho}(1))=0$.
By (\ref{eq:finite}) and hypothesis (\ref{pc:hyp2}), we conclude that
$\Hf(\GQ,A_{\rho}(1))=0$, as desired.
\end{proof}

\section{Local invariants, $\ell \neq p$} \label{s1}

\subsection{Characters} \label{s11}

Fix a prime $p$ and
let $\O \subseteq \C$ denote the ring of integers of a number
field $K$.  For each
prime $\lambda$ of $\O$ not dividing $p$, fix an isomorphism
$\ilm : \C \toi \Klmbar$ extending the inclusion $\O \inj \Olm$.  

Let $F$ be a finite extension of $\Qp$;
set $G_{F} = \Gal(\Fbar/F)$.  We say that a continuous character 
$\chi : F^{\times} \to \C^{\times}$ is of {\it Galois-type} with respect
to $\ilm$ if the character
$\ilm \circ \chi : F^{\times} \to \Klmbar^{\times}$
extends to a continuous character
$$\chilm : G_{F} \to \Klmbar^{\times}$$
via the dense injection $F^{\times} \inj G_{F}^{\ab}$ of local
class field theory.   We then write $\chiblm : G_{F} \to \Flbar^{\times}$ for
the reduction of $\chilm$.
We say that $\chi$ is
{\it arithmetic} if $\chi(F^{\times}) \subseteq \Qbar^{\times}$.

\begin{lemma} \label{lemma:trivial}
Let $\chi : F^{\times} \to \C^{\times}$ be an arithmetic character
of Galois-type
with respect to $\ilm$ for all $\lambda \nmid p$.  If
$\chiblm=1$ for infinitely many $\lambda$, then $\chi=1$.
\end{lemma}
\begin{proof}
Let $\bar{\lambda}$ be the prime of $\Qbar$ which is the kernel of the
composition
$$\O_{\Qbar} \overset{\ilm}{\too} \O_{\Klmbar} \to \Flbar.$$
If $\chiblm = 1$, then
$\chi(F^{\times}) \subseteq 1 + \bar{\lambda}$.  
If this hold for infinitely many
$\lambda$, then $\chi(F^{\times})-1$ lies in primes of $\Qbar$
of infinitely many distinct residue characteristics
and thus is trivial, as claimed.
\end{proof}

\subsection{Computation of local invariants} \label{s13}

Let $\pi$ be an irreducible admissible complex representation of
$\GL_{2}\Qp$.  (See \cite[Section 11]{DiamondIm} and the reference therein
for definitions, or 
see \cite{Buzzard} for a thoroughly enjoyable introduction.)
We say that $\pi$ is {\it arithmetic} if either:
\begin{enumerate}
\item $\pi$ is a subquotient of an induced representation
$\pi(\chi_{1},\chi_{2})$ with each $\chi_{i} : \Qp^{\times} \to \C^{\times}$
arithmetic; or
\item $\pi$ is the base change of an arithmetic 
character $\chi : F^{\times} \to
\C^{\times}$ for a quadratic extension $F/\Qp$; or
\item $\pi$ is extraordinary.
\end{enumerate}
We will recall the notion of {\it Langlands correspondence} (with
respect to $\ilm$) between
$\pi$ and a $\lambda$-adic representation $\rho : G_{p} \to
\GL_{2}\Klmbar$ in the course of the proof below.
We write $\nm : \Qp^{\times} \to
\C^{\times}$ for the norm character; it is arithmetic and of Galois-type
with respect to each $\ilm$,
with $\nm_{\lambda}$ equal to the $\lambda$-adic cyclotomic character 
$\eplm : G_{p} \to \Klm^{\times}$.

\begin{proposition} \label{prop:ss}
Let $\pi$ be an arithmetic irreducible admissible complex
representation of $\GL_{2}\Qp$.  Let $\{ \rholm : G_{p} \to \GL_{2}\Klmbar 
\}_{\lambda \nmid p}$ be a 
family of continuous representations such that $\pi$ and $\rholm$
are in Langlands correspondence with respect to $\ilm$ for
all $\lambda \nmid p$.  If $\pi$ is principal series or supercuspidal, then
$$H^{0}(G_{p},\eplmb \otimes \ad \rhoblm^{\sss}) = 0$$
for almost all $\lambda$.
\end{proposition}
\begin{proof}
For $\lambda$ of odd residue characteristic we have
$$\eplmb \otimes \ad \rhoblm^{\sss} \cong \eplmb \oplus \bigl(
\eplmb \otimes \adz \rhoblm^{\sss} \bigr).$$
The first summand is trivial if and only if $\lambda$ divides $p-1$,
so that we may restrict our attention to $\eplmb \otimes \adz \rhoblm^{\sss}$.

Assume first that $\pi = \pi(\chi_{1},\chi_{2})$ is arithmetic 
principal series.
In this case, for $\pi$ to be in Langlands correspondence with $\rholm$
means that the $\chi_{i} : \Qp^{\times} \to \C^{\times}$ are of Galois-type 
with respect to $\ilm$ and
$\rholm \cong \chi_{1,\lambda} \oplus \chi_{2,\lambda}$.
Thus
\begin{equation} \label{eq:ps}
\eplmb \otimes \adz \rhoblm^{\sss} \cong \eplmb \oplus 
\eplmb\chib^{\vphantom{-1}}_{1,\lambda}\chib_{2,\lambda}^{-1} \oplus
\eplmb\chib_{1,\lambda}^{-1}\chib^{\vphantom{-1}}_{2,\lambda}.
\end{equation}
We must show that each of these characters 
is non-trivial for almost all $\lambda$.
As above, this is clear for $\eplmb$.
If $\eplmb\chib^{\vphantom{-1}}_{1,\lambda}\chib_{2,\lambda}^{-1}$ is trivial
for infinitely many $\lambda$, then by Lemma~\ref{lemma:trivial} we must have
$\chi^{\vphantom{-1}}_{1}\chi_{2}^{-1} = \nm^{-1}$.  However,
$\pi$ would then be special or one dimensional, 
rather than principal series; thus this
can not occur.  The same argument deals with the
characters $\eplmb\chib_{1,\lambda}^{-1}\chib^{\vphantom{-1}}_{2,\lambda}$,
settling this case of the proposition.

Next assume that $\pi$ is supercuspidal but not extraordinary.  Then
$\pi$ is the base change of an arithmetic character $\chi : F^{\times} \to
\C^{\times}$ for a quadratic extension $F$ of $\Qp$.  The
Langlands correspondence in this case implies that $\chi$ is of Galois-type
with respect to each $\ilm$, and
$$\rholm \cong \Ind_{G_{F}}^{G_{p}} \chilm.$$
Let $\chi^{c}$ be the $\Gal(F/\Qp)$-conjugate character of $\chi$ and
let $\omega : F^{\times} \to \C^{\times}$
denote the character $\chi\cdot(\chi^{c})^{-1}$.  We have
$$\eplmb \otimes \adz \rhoblm^{\sss} \cong \eplmb\chi_{F} 
\oplus \bigl( \eplmb \otimes
\Ind_{G_{F}}^{G_{p}} \bar{\omega}_{\lambda} \bigr)$$
with $\chi_{F} : \Gal(F/\Qp) \to \{ \pm 1\}$ the non-trivial character
for $F/\Qp$.  The first summand is not a problem for $\ell > 3$.
The second summand
is irreducible if and only if $\bar{\omega}_{\lambda} \neq
\bar{\omega}^{\cc}_{\lambda}$.  In particular, if $\omega \neq \omega^{\cc}$,
then this case of
the proposition follows from Lemma~\ref{lemma:trivial}.
If instead $\omega = \omega^{\cc}$, then $\omega^{2} = 1$,
$\omega$ extends to a character
$\tilde{\omega} : \Qp^{\times} \to \C^{\times}$, and
$$\eplmb \otimes \adz \rhoblm^{\sss} \cong \eplmb\chi_{F} \oplus 
\eplmb\bar{\tilde{\omega}}_{\lambda} 
\oplus \eplmb\bar{\tilde{\omega}}^{-1}_{\lambda}.$$
Since $\omega^{2}=1$, we clearly have
$\tilde{\omega} \neq \nm^{\pm 1}$; this case of 
the proposition thus again follows from Lemma~\ref{lemma:trivial}.

Finally, if $\pi$ is extraordinary, then $p=2$, the image 
$\proj \rholm(I_{2})$ of inertia 
in $\PGL_{2}\Qlbar$ is isomorphic to $A_{4}$ or $S_{4}$,
and the composition
$$\proj \rholm(I_{2}) \inj \PGL_{2}\Qlbar \overset{\adz}{\too}
\GL_{3}\Qlbar$$
is an irreducible representation of $\proj \rholm(I_{2})$.
Since $\proj \rholm(I_{2})$ has order $12$ or $24$, it follows that
$\adz \rhoblm = \adz \rhoblm^{\sss}$ is an irreducible
$\Flbar$-representation of $I_{2}$ 
for $\lambda$ of residue characteristic at least $5$.  Thus already 
$H^{0}(I_{2},\eplmb \otimes \adz \rhoblm^{\sss}) = 0$
for such $\lambda$; the proposition follows.
\end{proof}

\begin{remark} \label{rmk:ss}
Note that if $\rhob : G_{p} \to \GL_{2}\Flbar$ is any reduction of $\rho$,
then
$$\dim_{\Flbar} H^{0}(G_{p},\eplmb \otimes \adz \rhob) \leq
\dim_{\Flbar} H^{0}(G_{p},\eplmb \otimes \adz \rhob^{\sss}).$$
\end{remark}

\begin{remark} \label{rmk:unr}
Suppose that $\pi = \pi(\chi_{1},\chi_{2})$ is an unramified principal series
representation (that is, $\chi_{i}(\Zp^{\times})=1$ for $i=1,2$).
Then if $\pi$ is in Langlands correspondence with $\rholm$, 
using (\ref{eq:ps}) one finds that:
\begin{multline*}
\quad H^{0}(G_{p},\eplmb \otimes \adz \rhoblm^{\sss}) \neq 0 \Rightarrow \\
p\bigl(\chi_{1}(p)+\chi_{2}(p)\bigr)^{2} \equiv 
(p+1)^{2}\chi_{1}(p)\chi_{2}(p) \pmod{\lambda}. \quad
\end{multline*}
In the case of a representation
$\rhoflm$ associated to a newform $f = \sum a_{n}q^{n}$ of 
weight $k$, level $N$, and character $\omega$ (see Section~\ref{s40}), this
translates to the condition:
$$H^{0}(G_{p},\eplmb \otimes \adz \rhob_{f,\lambda}^{\sss}) \neq 0
\,\,\Rightarrow\,\, 
a_{p}^{2} \equiv p^{k-2}(p+1)^{2}\omega(p) \pmod{\lambda}$$
for $p \nmid N\ell$.
\end{remark}

\begin{remark}
If $\pi$ is one dimensional or special, then
$$\dim_{\Flbar} H^{0}(G_{p},\eplmb \otimes \adz \rhob^{\sss}) = 1$$
for almost all $\lambda$.  For $\pi$ special, however, one might hope
to obtain an analogue of Proposition~\ref{prop:ss} by allowing $\rhob$ to
be non-semisimple.  Unfortunately, there is no way to formulate such a
result purely locally.  Instead, for the case of special local components
of Galois representations attached to
modular forms, in Section~\ref{s41}
we will use global considerations to rigidify $\rhob$; we will then be
able to obtain the required vanishing in the special case as well.
\end{remark}

\section{Local invariants, $\ell = p$} \label{s2}

\subsection{Filtered Dieudonn\'e modules} \label{s21}

We briefly review the theory of Fontaine--Laffaille \cite{FontaineLaffaille}.
Let $K$ be a finite extension of $\Ql$ with ring of integers $\O$.
Let $v : K^{\times} \to \Q$ be the valuation on $K$, normalized so that
$v(\ell) = 1$; let $\lambda$ be a uniformizer of $\O$ and 
let $e = v(\lambda)^{-1}$ denote the absolute ramification degree of $K$.

\begin{definition} \label{def:fdm}
A {\it filtered Dieudonn\'e $\O$-module (over $\Zl$)} is an
$\O$-module $D$ of finite type,
endowed with a decreasing filtration
$(D^{i})_{i \in \Z}$ by $\O$-module direct summands
and a family $(f_{i} : D^{i} \to D)_{i \in \Z}$ of $\O$-linear maps
satisfying:
\begin{enumerate}
\item \label{dm2} $D^{i}=D$ (resp.\ $D^{i}=0$) 
for $i \ll 0$ (resp.\ $i \gg 0$);
\item \label{dm3} $f_{i}|_{D^{i+1}} = \ell \cdot f_{i+1}$ for all $i$;
\item \label{dm4} $D = \sum_{i \in \Z} f_{i}(D^{i})$.
\end{enumerate}
For $a < b$, let $\MF^{a,b}(\O)$ 
denote the category of filtered Dieudonn\'e $\O$-modules $D$
satisfying $D^{a}=D$ and $D^{b}=0$.
\end{definition}

We now recall the relation between filtered Dieudonn\'e $\O$-modules and
Galois representations.  Let $V$ be a finite dimensional $K$-vector space
with a continuous $K$-linear action of $G_{\ell}$.  Define
a $K$-vector space
$$D_{\crys}(V) := (B_{\crys} \otimes_{\Ql} V)^{G_{\ell}}$$
with $B_{\crys}$ the crystalline period ring of Fontaine; $D_{\crys}(V)$
inherits a decreasing filtration $(D^{i}_{\crys}(V))_{i \in \Z}$
from the filtration on $B_{\crys}$.
We say that $V$ is {\it crystalline} if $D_{\crys}(V)$ and $V$ are
$K$-vector spaces of the same dimension.

For $a < b$,
let $\G^{a,b}(\O)$ denote the category of finite type $\O$-module subquotients
of crystalline $K$-representations $V$ with 
$D_{\crys}^{a}(V) = D_{\crys}(V)$ and
$D_{\crys}^{b}(V)=0$.  Fontaine--Laffaille define a functor
$$\U : \MF^{a,a+\ell}(\O) \to \G^{a,a+\ell}(\O)$$
which is equivalent to the identity functor on the underlying $\O$-modules and 
which induces an equivalence of categories between
$\MF^{a,a+\ell-1}(\O)$ and $\G^{a,a+\ell-1}(\O)$.
(See \cite[Section 1.1]{DFG} for more details; 
note that we are using Tate twists as in \cite[Section 4]{BlochKato}
to extend $\U$ to the case $a \neq 0$.)

\begin{example} \label{ex:char}
Let $\omega : G_{\ell} \to \O^{\times}$ be 
an unramified character of finite order and
let $\O(\omega)$ denote a free $\O$-module of rank $1$ with 
$G_{\ell}$-action via
$\omega$.  Then $\O(\omega) \in \G^{0,1}(\O)$, so that
there is $D_{\omega} \in \MF^{0,1}(\O)$ such that $\U(D_{\omega})\cong
\O(\omega)$.
This $D_{\omega}$ is a free $\O$-module of rank one with
$D_{\omega} = D_{\omega}^{0}$ and
$f_{0}$ multiplication by $\omega^{-1}(\ell)$.
\end{example}

\begin{example} \label{ex:mf}
Let $f = \sum a_{n}q^{n}$ be a newform of weight $k \geq 2$, level $N$, 
and character $\omega$.  Assuming that $K$ contains some completion of the
number field generated by the $a_{n}$, there is a Galois representation
$\rho_{f} : \GQ \to \GL_{2}K$ associated to $f$ as in Section~\ref{s40}.
Fix an embedding $G_{\ell} \inj \GQ$ and
let $V_{f}$ be a two dimensional $K$-vector space on which $G_{\ell}$ acts
via $\rho_{f}|_{G_{\ell}}$.  

Fix a $G_{\ell}$-stable $\O$-lattice $T_{f}
\subseteq V_{f}$.  If $\ell \nmid N$, then $V_{f}$ is crystalline and
$T_{f} \in \G^{0,k}(\O)$.
Thus for $\ell > k$ there exists $D_{f} \in \MF^{0,k}(\O)$ with
$\U(D_{f})\cong T_{f}$.  
Using \cite[Theorem 4.3]{BlochKato} and standard properties of modular
representations, 
one obtains the following description of $D_{f}$.  The filtration satisfies:
$$\rank_{\O} D_{f}^{i} = \begin{cases} 2 & i \leq 0; \\ 
1 & 1 \leq i \leq k-1; \\
0 & k \leq i. \end{cases}$$
Choose an $\O$-basis $x,y$ of $D_{f}$ with $x$ an $\O$-generator
of $D_{f}^{1}$.  Let $\alpha,\beta,\gamma,\delta \in \O$ be
such that
$$f_{0}x = \alpha x + \beta y; \qquad f_{0}y = \gamma x + \delta y.$$
Then $\alpha + \delta = a_{\ell}$ and $\alpha\delta - \beta\gamma = 
\ell^{k-1}\omega(\ell)$.
By (\ref{dm3}) of Definition~\ref{def:fdm}
we have $v(\alpha),v(\beta) \geq k-1$.
\end{example}

\subsection{Computation of local invariants} \label{s22}

Let $f$ and $T_{f}$
be as in Example~\ref{ex:mf}.  Let
$$\rhob_{f} : G_{\ell} \to \GL_{2} (\O/\lambda)$$
be the Galois representation on $T_{f}/\lambda T_{f}$.

\begin{proposition} \label{prop:adjdm}
Assume $\ell \nmid N$ and $\ell > 2k$.  Then
$$H^{0}(G_{\ell},\eplmb \otimes\ad \rhob_{f}) = 0$$
unless $k=2$ and $a_{\ell}^{2} \equiv \omega(\ell) \pmod{\lambda}$.
\end{proposition}
\begin{proof}
Since $\det \rho_{f} = \eplm^{1-k}\omega^{-1}$, we have
\begin{equation} \label{eq:decomp}
\ad T_{f}(1) \cong \bigl(T_{f} \otimes_{\O} T_{f} \otimes_{\O} 
\O(\omega)\bigr)(k)
\end{equation}
(where $(-)(k)$ denotes the $k$-fold Tate twist).
It follows that $\ad T_{f}(1) \in \G^{-k,k-1}(\O)$.  Since $\ell > 2k$, there
thus exists $D \in \MF^{-k,k-1}(\O)$ with $\U(D) \cong \ad T_{f}(1)$.
In fact, by (\ref{eq:decomp}) and 
\cite[Proposition 1.7]{FontaineMessing} we can take
$$D = \bigl(D_{f} \otimes_{\O} D_{f} \otimes_{\O}
D_{\omega}\bigr)(k)$$
with notation as in Examples~\ref{ex:char} and~\ref{ex:mf}.
Further, since $\ad T_{f}(1)/\lambda$ is a realization of
$\eplmb \otimes \ad \rhob_{f}$, 
by \cite[Lemma 4.5]{BlochKato} we have
\begin{equation} \label{eq:coh}
H^{0}(G_{\ell},\eplmb \otimes \ad \rhob_{f}) \cong
\ker (1 - f_{0} : D^{0}/\lambda D^{0} \to D/\lambda D).
\end{equation}
Thus to prove the proposition it suffices to compute the latter group.

By the definition of Tate twists
and tensor products of filtered Dieudonn\'e $\O$-modules, we have
\begin{align*}
D^{0} &= (D_{f} \otimes D_{f} \otimes D_{\omega})^{k} \\
&= \sum_{i_{1}+i_{2}+i_{3}=k} D_{f}^{i_{1}} \otimes D_{f}^{i_{2}} \otimes
D_{\omega}^{i_{3}} \\
&= D_{f}^{1} \otimes D_{f}^{k-1} \otimes D_{\omega}^{0} \\
&= \O \cdot (x \otimes x \otimes w)
\end{align*}
where $x$ is as in Example~\ref{ex:mf} and $w$ is an $\O$-generator of
$D_{\omega}$.  Using (\ref{dm3}) of Definition~\ref{def:fdm}, we compute:
\begin{align}
f_{0}(x \otimes x \otimes w) &:= 
f_{1}x \otimes f_{k-1}x \otimes f_{0}w \notag \\
&= \frac{\omega^{-1}(\ell)}{\ell^{k}}
(\alpha x + \beta y) \otimes (\alpha x + \beta y) \otimes w. \label{eq:ff}
\end{align}

Suppose now that (\ref{eq:coh}) is non-zero.  Then by (\ref{eq:ff}) we must
have
\begin{gather}
\omega^{-1}(\ell)\alpha^{2} \equiv \ell^{k} \pmod{\lambda^{ek+1}}; 
\label{eq:big}\\
\alpha\beta \equiv 0 \pmod{\lambda^{ek+1}}. \label{eq:val}
\end{gather}
In particular, we must have $v(\alpha) = k/2$.  Since also
$v(\alpha) \geq k-1$, this implies that $k=2$ and $v(\alpha)=1$.
Thus (\ref{eq:val}) implies that $v(\beta) > 1$.  As
$\alpha\delta - \beta\gamma = \ell\omega(\ell)$,
we conclude that $v(\delta)=0$ and
$$\alpha\delta \equiv \ell\omega(\ell) \pmod{\lambda^{e+1}}.$$
Using this and (\ref{eq:big}) one deduces easily that
$$a_{\ell}^{2} = (\alpha + \delta)^{2} \equiv \omega(\ell) \pmod{\lambda}$$
as claimed.
\end{proof}

\section{Unobstructed deformation problems} \label{s4}

\subsection{Modular representations} \label{s40}

Let $f = \sum a_{n}q^{n}$ be a newform of weight $k \geq 2$,
level $N$, and character $\omega$. Fix a finite set of primes $S$ containing
all primes dividing $N$.  Let $K$ be the number field
generated by the $a_{n}$.  Then for any prime $\lambda$ of $K$
(say with residue field $\klm$ of characteristic $\ell$)
there is a continuous Galois representation
$$\rhoflm : \GQSl \to \GL_{2}\Klm,$$
unramified with $\trace \rhoflm(\Frob_{p}) = a_{p}$
for $p \nmid N\ell$; the determinant of $\rhoflm$ is
$\eplm^{1-k}\omega^{-1}$.
(We find it more convenient to work with this geometric normalization of
$\rhoflm$, which is dual to the more common arithmetic normalization.)
Let
$$\rhobflm : \GQSl \to \GL_{2}\klm.$$
be the semisimple reduction of $\rhoflm$.
Then $\rhobflm$ is absolutely irreducible for almost all
primes $\lambda$ by \cite[Lemma 7.13]{DFG}.

Let $\pi$ be the automorphic representation corresponding to $f$ as in
\cite[Section 11.1]{DiamondIm}
and let $\pi = \otimes' \pi_{p}$ be the decomposition of $\pi$ into
admissible complex representations $\pi_{p}$ of $\GL_{2} \Qp$.
Fix isomorphisms $\ilm : \C \toi \Klmbar$ (extending the inclusion
$K \inj \Klm$) for all $\lambda$.
By \cite[Th\'eor\`eme B]{Carayol} and
\cite[Proposition 9.3]{JacquetLanglands}
each
$\pi_{p}$ is arithmetic and infinite dimensional and is
in Langlands correspondence (with respect to $\ilm$) with
$\rho_{f,\lambda}|_{G_{p}}$ for each $\lambda \nmid p$.

\subsection{Special primes} \label{s41}

We continue with the notation of the previous section.
Let $p$ be a prime such that $\pi_{p}$ is special; that is, $\pi$ is
the unique infinite dimensional subquotient of $\pi(\chi\nm,\chi)$ for
an arithmetic character $\chi : \Qp^{\times} \to \C^{\times}$.  Then
$\chi$ is of Galois-type
with respect to $\ilm$ for each $\lambda \nmid p$, and
$$\rhoflm|_{G_{p}} \cong \left(\begin{array}{cc} \eplm\chilm & * \\
0 & \chilm \end{array}\right)$$
with the upper right-hand corner non-zero and ramified.

\begin{lemma} \label{lemma:woops}
If $p^{2} \not\equiv 1 \pmod{\lambda}$, then
$$\rhobflm|_{G_{p}} \otimes \klmbar \cong
\left(\begin{array}{cc} \eplmb\chiblm & * \\
0 & \chiblm \end{array}\right).$$
\end{lemma}
\begin{proof}
Since the semisimplfication of $\rhobflm|_{G_{p}} \otimes \klmbar$ is 
$\eplmb\chiblm \oplus \chiblm$, the
only way the lemma can fail is if
$$\rhobflm|_{G_{p}} \otimes \klmbar \cong
\left(\begin{array}{cc} \chiblm & \nu \\
0 & \eplmb\chiblm \end{array}\right)$$
with $\nu$ non-trivial.  One checks directly that $\eplmb^{-1}\chiblm^{-1}\nu$ 
is naturally
an element of $H^{1}(G_{p},\klmbar(-1))$; since this cohomology group is 
trivial
unless $p^{2} \equiv 1 \pmod{\lambda}$, the lemma follows.
\end{proof}

Using Lemma~\ref{lemma:woops},
the proof of the next lemma is a straightforward matrix computation; we omit
the details.

\begin{lemma} \label{lemma:spec}
Let $p$ be a prime for which $\pi_{p}$ is special.  
Let $\lambda$ be a prime not dividing $2p(p^{2}-1)$.  Then
$$H^{0}(G_{p},\eplmb \otimes \adz \rhobflm) \neq 0$$
if and only if $\rhobflm|_{G_{p}} \otimes \klmbar$ is semisimple.
\end{lemma}

The level-lowering approach of the next proposition was suggested to the
author by Ken Ribet.

\begin{proposition} \label{prop:special}
Let $p$ be a prime for which $\pi_{p}$ is special.  Then
$$H^{0}(G_{p},\eplmb \otimes \ad \rhobflm) = 0$$
for almost all primes $\lambda$.
\end{proposition}
\begin{proof}
As always it suffices to prove the proposition for $\eplmb \otimes
\adz \rhobflm$.  Let $\chi$ be as above.
Since the automorphic representation
$\pi$ has central quasi-character $\nm^{1-k}\omega^{-1}$, we must have
$\chi^{2} = \nm^{-k}\omega_{p}^{-1}$, 
where $\omega_{p} : \Qp^{\times} \to
\C^{\times}$ is the $p$-component of $\omega$.
In particular, $\chi_{p}' := \chi^{-1}\nm^{-k/2}$ has finite order.
Extend $\chi'_{p}$ to a Dirichlet character $\chi'$ and
let $f'$ denote the newform of weight $k$ and some level $M$
associated to the eigenform $f \otimes \chi'$.
The $p$-component $\pi_{p}'$ of the automorphic representation
associated to $f'$ is a subquotient of $\pi(\chi'_{p}\chi\nm,\chi'_{p}\chi)$;
since $\chi'_{p}\chi$ is unramified at $p$ by construction, we conclude that
$p$ divides $M$ exactly once.

Suppose now that $\lambda \nmid 2p(p^{2}-1)$ is such $\rhobflm$ is irreducible
and
\begin{equation} \label{eq:last}
H^{0}(G_{p},\eplmb \otimes \adz \rhobflm) \neq 0.
\end{equation}
Then $\rhobflm|_{G_{p}} \otimes \klmbar$ 
is semisimple by Lemma~\ref{lemma:spec}, so that
$$\rhobflm|_{G_{p}} \otimes \klmbar \cong \eplmb\chiblm \oplus \chiblm$$
by Lemma~\ref{lemma:woops}.  Thus
$$\rhob_{f',\lambda}|_{G_{p}} \cong \rhobflm|_{G_{p}} \otimes \chi'_{p}
\cong \eplmb^{1-\frac{k}{2}} \oplus \eplmb^{-\frac{k}{2}}.$$
In particular $\rhob_{f',\lambda}$ is unramified at $p$.
By \cite[Theorem 1.1]{Diamond} applied to $f'$
there is then a newform $f''$ of weight $k$ and level dividing
$M/p$ with $f'$ congruent to 
$f''$ modulo some prime of $\Qbar$ above $\lambda$.  However, there are only
finitely many newforms $f''$ of weight $k$ and level
dividing $M/p$, and for each such $f''$ there are
only finitely many $\lambda$ with $f'$ congruent to $f''$ modulo some prime
of $\Qbar$ above $\lambda$.  
Thus (\ref{eq:last}) can only hold for finitely many $\lambda$;
the proposition follows.
\end{proof}

\subsection{Proof of Theorem~\ref{mt1}} \label{s42}

Let $f$ and $S$ be as above.
For any $\lambda$ at which $\rhobflm$ is absolutely irreducible, we
write $R^{S}_{f,\lambda}$ for the universal deformation ring for
the deformation problem $D^{S \cup \{\ell\}}_{\rhobflm}$
as in Section~\ref{s31}.

\begin{theorem} \label{thm:main1}
Suppose $k > 2$.  Then
$$R^{S}_{f,\lambda} \cong W(\klm)[[T_{1},T_{2},T_{3}]]$$
for almost all primes $\lambda$.
\end{theorem}
\begin{proof}
We apply Proposition~\ref{prop:criterion} to the lifting $\rhoflm$ of
$\rhobflm$.  Hypothesis (1) holds
for almost all $\lambda$ by Proposition~\ref{prop:ss} and 
Remark~\ref{rmk:ss} (for $p \in S$
for which $\pi_{p}$ is not special), Proposition~\ref{prop:special} (for
$p \in S$ at which $\pi_{p}$ is special), and Proposition~\ref{prop:adjdm}
(for $p = \ell$).  Hypotheses (2) and (3) hold for almost all $\lambda$
by \cite[Theorems 8.2 and 7.15]{DFG} and the fact that
(in the notation of \cite{DFG})
$$\Hf(\GQ,A_{\rhoflm}) \subseteq H^{1}_{\Sigma}(\GQ,A_{\rhoflm})$$
for any finite set of primes $\Sigma$.
Thus the deformation problem $D^{S \cup \{\ell\}}_{\rhobflm}$
is unobstructed for almost all $\lambda$; Proposition~\ref{prop:defthy}
completes the proof.
\end{proof}

\begin{theorem} \label{thm:main2}
Suppose $k = 2$.  Then
\begin{equation} \label{eq:almost}
R^{S}_{f,\lambda} \cong W(\klm)[[T_{1},T_{2},T_{3}]]
\end{equation}
for almost all primes $\lambda$ dividing $\ell$ with
$a_{\ell}^{2} \not\equiv \omega(\ell) \pmod{\lambda}.$
In particular, (\ref{eq:almost}) holds for a set of $\lambda$ of
density one.
\end{theorem}
\begin{proof}
The proof of the first statement 
is the same as the proof of Theorem~\ref{thm:main1}, taking into
account the modifications for $k=2$ in Proposition~\ref{prop:adjdm}.
The density statement follows from
\cite[Theorem 20]{Serre} and the Ramanujan--Petersson conjecture
(proven by Deligne).
\end{proof}

Note that by \cite[Section 1.3]{Mazur}, the analogue of these
results for the arithmetic normalization of $\rhobflm$ hold as well.

\subsection{Modular forms of level one} \label{s43}

Let $f$ be the unique normalized cusp form of level $1$, weight 
$k = 12,16,18,20,22 \text{~or~}26$,
and trivial character.  Then $f$ has rational
coefficients, so that for every prime $\ell$ we obtain a
representation
$$\rhobfl : \GQl \to \GL_{2}\Fl.$$
We recall the set of primes $\ell$ for which $\rhobfl$ is not absolutely
irreducible (see \cite{Serre2}):

\vspace{0.2cm}
\begin{center}
\begin{tabular}{c|l||c|l}
$k$ & $\ell$ & $k$ & $\ell$ \\ \hline
12 & 2,3,5,7,691 & 20 & 2,3,5,7,11,13,283,617 \\
16 & 2,3,5,7,11,3617 & 22 & 2,3,5,7,13,17,131,593 \\
18 & 2,3,5,7,11,13,43867 & 26 & 2,3,5,7,11,17,19,657931 
\end{tabular}
\end{center}
\vspace{0.2cm}

\begin{theorem} \label{thm:one}
Let $f$ be as above and let $\ell > k+1$ be a prime for which $\rhobfl$ is
absolutely irreducible.  Then
\begin{equation} \label{eq:one}
R_{f,\ell}^{\emptyset} \cong \Zl[[T_{1},T_{2},T_{3}]].
\end{equation}
\end{theorem}

We leave it to the reader to use Remark~\ref{rmk:unr} to derive
the generalization of Theorem~\ref{thm:one} to the case $S \neq \emptyset$.

\begin{proof}
Fix $\ell > 2k$.  Then by \cite[Corollary 7.2 and Section 7.4]{DFG} we have
\begin{equation} \label{eq:cong}
\# \Hf(\GQ,A_{\rhofl}) =
\# H^{1}_{\emptyset}(\GQ,A_{\rhofl}) = \# \Zl/\eta_{f}^{\emptyset}\Zl
\end{equation}
where $\eta_{f}^{\emptyset} \subseteq \Z$ is the congruence ideal of
\cite[Section 6.4]{DFG}.  By definition, $\eta_{f}^{\emptyset}$ is generated
by the element $d(f)$ of \cite[Theorem 1]{Hida}.
Using \cite[Theorem 1 and Theorem 2]{Hida} and the fact that $f$
is the unique normalized cusp form of its weight, level, and character, 
we see that $d(f)$ is divisible only by primes $\leq k - 2$.  In particular,
$d(f)$ is not divisible by $\ell$, so that (\ref{eq:cong}) implies that
\begin{equation} \label{eq:van}
\Hf(\GQ,A_{\rhofl}) = 0.
\end{equation}

Since $S = \emptyset$, Propositions~\ref{prop:criterion} and~\ref{prop:adjdm}
now imply that (\ref{eq:one}) holds for
$\ell > 2k$.  The fact that it holds for $k + 1 < \ell < 2k$
follows on checking the criterion of \cite[Section 7]{Mazur3}
for each such $\ell$.
\end{proof}


\providecommand{\bysame}{\leavevmode\hbox to3em{\hrulefill}\thinspace}

\end{document}